\documentclass[12pt,a4paper]{article}
\usepackage[T1]{fontenc}
\usepackage{mathptmx}
\usepackage{courier}
\usepackage[dvips]{graphicx}
\usepackage{psfrag}
\usepackage{amsmath,amssymb}        
\usepackage{latexsym}
\usepackage{url}
\usepackage[light]{kpfonts}

\def\ds{\displaystyle}

\title{An Iterative Transformation Method\\ for a Similarity Boundary Layer Model}
\author{Riccardo Fazio\\
Department of Mathematics, Computer Science,\\ Physical Sciences and Earth Sciences,\\
University of Messina \\
Viale F. Stagno D'Alcontres, 31 \\
98166 Messina, Italy \\
E-mails: \url{rfazio@unime.it} \\
Home-page: \url {http://mat521.unime.it/fazio}} 
\date{\today}
\begin{document}
\maketitle
\begin{abstract}
In this paper, within scaling invariance theory, we define and apply to the numerical solution of a similarity boundary layer model an iterative transformation method.
The boundary value problem to be solved depends on a parameter and is defined on a semi-infinite interval.
Using our transformation method we are able to solve the problem in point for a large range of the parameter involved.
As far as the accuracy of our numerical method is concerned, for two specific values of the involved parameter, we are able to compare favorably the obtained numerical result for the so-called missing initial condition to the exact solution reported by Crane [Z. Angew. Math. Phys., 21:645-647, 1970].
\end{abstract}
\bigskip

\bigskip


\section{Introduction}
The first definition and application of a transformation method dates back to 1912 when T\"opfer \cite{Topfer:1912:BAB} published a paper where he reduced the solution of the Blasius problem to the solution of a related initial value problem (IVP). 
From that time non-Iterative Transformation Methods (ITM) have been applied to several problems of practical interest within the applied sciences. 
First of all, a non-ITM was applied to the Blasius equation with slip boundary condition, arising within the study of gas and liquid flows at the micro-scale regime \cite{Gad-el-Hak:1999:FMM,Martin:2001:BBL}, see \cite{Fazio:2009:NTM}.
A non-ITM was applied also to the Blasius equation with moving wall considered by Ishak et al. \cite{Ishak:2007:BLM} or surface gasification studied by Emmons \cite{Emmons:1956:FCL} and recently by Lu and Law \cite{Lu:2014:ISB} or slip boundary conditions investigated by Gad-el-Hak \cite{Gad-el-Hak:1999:FMM} or Martin and Boyd \cite{Martin:2001:BBL}, see Fazio \cite{Fazio:2016:NIT} for details.
In particular, within these applications, we found a way to solve non-iteratively the Sakiadis problem \cite{Sakiadis:1961:BLBa,Sakiadis:1961:BLBb}.
The application of a non-ITM to an extended Blasius problem has been the subject of a recent manuscript \cite{Fazio:2020:NIT}.
A recent review dealing with the non-ITM and its applications can be be found, by the interested reader, in \cite{Fazio:2019:NIT}.

A first extension T\"opfer's to classes of problems in boundary layer theory involving one or more physical parameters
was considered by Na \cite{Na:1970:IVM}, see also NA \cite[Chapters 8-9]{Na:1979:CME} for an extensive survey on this subject.

Finally, an iterative extension of T\"opfer's algorithm has been introduced, for the numerical solution of free BVPs, by Fazio \cite{Fazio:1990:SNA}. 
This iterative extension has been applied to several problems of interest: free boundary problems \cite{Fazio:1990:SNA,Fazio:1997:NTE,Fazio:1998:SAN},
a moving boundary hyperbolic problem \cite{Fazio:1992:MBH}, Homann and Hiemenz problems governed by the Falkner-Skan equation and a mathematical model describing the study of the flow of an incompressible fluid around a slender parabola of revolution \cite{Fazio:1994:FSE,Fazio:1996:NAN},
one-dimensional parabolic moving boundary problems \cite{Fazio:2001:ITM}, two variants of the Blasius problem \cite{Fazio:2009:NTM}, namely: a boundary layer problem over moving surfaces, studied first by Klemp and Acrivos \cite{Klemp:1972:MBL}, and a boundary layer problem with slip boundary condition, that has found application in the study of gas and liquid flows at the micro-scale regime \cite{Gad-el-Hak:1999:FMM,Martin:2001:BBL}, parabolic problems on unbounded domains \cite{Fazio:2010:MBF} and, recently, see \cite{Fazio:2015:ITM}, a further variant of the Blasius problem in boundary layer theory: the so-called Sakiadis problem \cite{Sakiadis:1961:BLBa,Sakiadis:1961:BLBb}.
Moreover, this iterative extension can be used to investigate the existence and uniqueness question for different class of problems, as shown for free BVPs in \cite{Fazio:1997:NTE}, and for problems in boundary layer theory in \cite{Fazio:2020:EUB}.
A recent review dealing with, the derivation and applications of the ITM can be be found, by the interested reader, in \cite{Fazio:2019:ITM}.
A unifying framework, providing proof that the non-ITM is a special instance of the ITM and consequently can be derived from it, has been the argument of the paper \cite{Fazio:2020:SIT}.

\section{Mathematical Model}
In this section we report the boundary value problem (BVP) that we would like to solve.
This BVP is defined by
\begin{align}\label{eq:goveq} 
& {\displaystyle \frac{d^{3}f}{d\eta^3}} + f
{\displaystyle \frac{d^{2}f}{d\eta^2}} - \beta \left(\frac{df}{d\eta}\right)^2 = 0 \nonumber \\[-1ex]
&\\[-1ex]
& f(0) = 0 \ , \qquad {\displaystyle \frac{df}{d\eta}}(0) = 1 \ , \qquad {\displaystyle \frac{df}{d\eta}}(\eta) \rightarrow 0 \quad \mbox{as}
\quad \eta \rightarrow \infty \ , \nonumber 
\end{align}
where $\beta$ is a given value verifying the limitations: $-2 < \beta \le 2$ for the free convection flow over a vertical semi-infinite flat plate, see Liao and Pop \cite{Liao:2004:EAS}, and 
$ 2 < \beta < +\infty $ in the case of boundary layer flows over a stretching wall, see Crane \cite{Crane:1970:FPS}, Banks \cite{Banks:1983:SSB}, Magyari and Keller \cite{Magyari:2000:ESS}, Brown and Stewartson \cite{Brown:1965:SSB} and Stuart \cite{Stuart:1966:DBL}. 
Banks \cite{Banks:1983:SSB} showed that there is no solution for $\beta = -2$, and this was rigorously proved by Ingham and Brown \cite{Ingham:1986:FPS}.
Moreover, Ingham and Brown \cite{Ingham:1986:FPS} reported that they found, numerically, the existence of dual solutions, for the same value of $\beta > 0$.
The second branch of solution having negative values of the velocity $\frac{df}{d\eta}(\eta)$ in some regions and this is not realistic from a physical point of view.

\section{The ITM} 
In order to define our ITM we need to embedd the original BVP (\ref{eq:goveq}) into a class of BVPs depending on a numerical parameter $h$
\begin{align}\label{eq:EXgoveq} 
& {\displaystyle \frac{d^{3}f}{d\eta^3}} + f
{\displaystyle \frac{d^{2}f}{d\eta^2}} - \beta \left(\frac{df}{d\eta}\right)^2 = 0 \nonumber \\[-1ex]
&\\[-1ex]
& f(0) = 0 \ , \qquad {\displaystyle \frac{df}{d\eta}}(0) = h^{(1-\delta)/\sigma} \ , \qquad {\displaystyle \frac{df}{d\eta}}(\eta) \rightarrow 1-h^{(1-\delta)/\sigma} \quad \mbox{as}
\quad \eta \rightarrow \infty \ , \nonumber 
\end{align}
and consider the invariance of the governing differential equation and of the prescribed initial conditions with respect to the extended scaling group of point transformations
\begin{equation}\label{eq:scaling2}
f^* = \lambda f \ , \qquad \eta^* = \lambda^{\delta} \eta \ , \qquad h^* = \lambda^{\sigma} h \ .   
\end{equation}
It is easily seen, that the governing differential equation and the prescribed initial conditions are invariant on condition that $\delta = -1$.
Moreover, the introduced scaling group involves the scaling of the fictitious numerical parameter $h$, that have been used to force the initial conditions to be invariant and the asymptotic boundary condition to be not invariant. 
Now, we can integrate the governing equation in (\ref{eq:goveq}) written in the star variables on $[0, \eta^*_\infty]$, where $\eta^*_\infty$ is a suitable truncated boundary, with initial conditions
\begin{equation}\label{eq:ICs2}
f^*(0) = 0 \ , \qquad \frac{df^*}{d\eta^*}(0) = h^{2/\sigma}\ , \quad {\displaystyle \frac{d^2f^*}{d\eta^{*2}}(0)} = -1 \ ,
\end{equation}
in order to compute an approximation $\frac{df^*}{d\eta^*}(\eta^*_{\infty})$ for $\frac{df^*}{d\eta^*}(\infty)$ and the corresponding value of $\lambda$ by the equation
\begin{equation}\label{eq:lambda}
\lambda = \left[\frac{df^*}{d\eta^*}(\eta_\infty^*) + h^{*2/\sigma}\right]^{1/2} \ .
\end{equation}
Once the value of $\lambda$ has been computed by equation (\ref{eq:lambda}), we can find the missed initial condition by the relation
\begin{equation}\label{eq:MIC}
\frac{d^2f}{d\eta^{2}}(0) =  \lambda^{-3} \frac{d^2f^*}{d\eta^{*2}}(0) \ .
\end{equation}
In the ITM we proceed as follows: we set the values of $\sigma$ and $\eta_{\infty}^*$ and integrate the IVP on $[0, \eta_{\infty}^*]$.
Naturally, choosing $h^*$ arbitrarily we do not obtain the value $h = 1$, however, we can apply a root-finder method, like bisection, secant, regula-falsi, Newton or quasi-Newton root-finder, because the required value of $h$ can be considered as the root of the implicit defined, transformation function
\begin{equation}\label{eq:Tfunction}
\Gamma(h^*) = \lambda^{-\sigma} h -1  \ .   
\end{equation}
Of course, any positive value of $\sigma$ can be chosen, and in the following, for the sake of simplicity, we set $\sigma = 4$.

\section{Numerical results}
In this section, we report the numerical results obtained by the ITM.
As a root-finder we apply Newton method, with a termination criterion given by $|\Gamma(h^*)| < Tol$ with $Tol = 10^{-9}$.
Moreover, some preliminary numerical experiments allowed us to set, the truncated boundary value, $\eta^*_{\infty} = 5$.
However, the application of Newton's root-finder requires a more complex treatment involving a system of six differential equations. 
In fact, in order to apply the Newton's root-finder, at each iteration, we have to compute the derivative
with respect to $h^*$ of the transformation function $\Gamma$.
In our case, replacing equation (\ref{eq:lambda}) into (\ref{eq:Tfunction}), the transformation function is given by 
\begin{equation}\label{eq:Gamma}
\Gamma(h^*) = \left[u_2^*(\eta_\infty^*)+h^{*1/2}\right]^{-2} h^* -1 \ ,
\end{equation}
and its first derivative can be easily computed as
\begin{equation}\label{eq:dGamma}
\frac{d\Gamma}{dh^*}(h^*) = \left[u_2^*(\eta_\infty^*)+h^{*1/2}\right]^{-2}\left\{1-2 \left[u_5^*(\eta_\infty^*)+\frac{1}{2}h^{*-1/2}\right]\left[u_2^*(\eta_\infty^*)+h^{*1/2}\right]^{-1} h^*\right\}  \ .
\end{equation}

Let us introduce the auxiliary variables $u_j(\eta)$ for $j = 1, 2, \dots , 6$ defined by
\begin{equation}\label{eq:var}
u_1 = f \ , \quad u_2 = \frac{df}{d\eta} \ , \quad u_3 = \frac{d^2f}{d\eta^2} \ , \quad u_4 = \frac{\partial u_1}{\partial h^*} \ , \quad u_5 = \frac{\partial u_2}{\partial h^*} \ , \quad u_6 = \frac{\partial u_3}{\partial h^*} \ .
\end{equation}
Now, the related IVP is given by
\begin{align}\label{eq:IVP6}
& \frac{du_1^*}{d\eta^*} = u_2^* \ , \nonumber \\
& \frac{du_2^*}{d\eta^*} = u_3^* \ , \nonumber \\
& \frac{du_3^*}{d\eta^*} = \beta u_2^* u_2^* -u_1^* u_3^* \ , \nonumber \\[-1.5ex]
&\\[-1.5ex]
& \frac{du_4^*}{d\eta^*} = u_5^* \ , \nonumber \\
& \frac{du_5^*}{d\eta^*} = u_6^* \ , \nonumber \\
& \frac{du_6^*}{d\eta^*} = 2 \beta u_2^* u_5^* -u_4^* u_3^* - u_1^* u_6^* \ , \nonumber 
\end{align}
$u_1^*(0) = 0 \ , \ u_2^*(0) = h^{*1/2} \ , \ u_3^*(0) = -1 \ , \ 
u_4^*(0) = 0 \ , \ u_5^*(0) = \frac{1}{2} h^{*-1/2} \ , \ u_6^*(0) = 0$.

In table \ref{tab:ITM:ITERA} we report a sample iteration of the Newton's method.
\begin{table}[!hbt]
\caption{Newton's method iterations for $h_0^* = 1.75$. Here the $D$ notation stands for a double precision arithmetic.}
\vspace{.5cm}
\renewcommand\arraystretch{1.3}
	\centering
		\begin{tabular}{clr@{.}lr@{.}l}
\hline 
iteration & {$h^*$} & \multicolumn{2}{c}%
{$\lambda$} & \multicolumn{2}{c}%
{$\Gamma(h^*)$} \\[1.2ex]
\hline
0 & 1.75     &  1  & 108575 &    0 & 158719908 \\ 
1 & 1.837475 &  1  & 157106 &    0 & 025010170 \\
2 & 1.856888 &  1  & 167093 &    0 & $84D-04$ \\
3 & 1.857586 &  1  & 167447 &    1 & $01D-06$ \\
4 & 1.857587 &  1  & 167447 &    1 & $59D-12$ \\
\hline			
		\end{tabular}
	\label{tab:ITM:ITERA}
\end{table}
As it is easily seen from this table, as a result of the second order convergence of Newton's root-finder we get the chosen termination criterion verified after only four iterates.
Moreover, we found that the value of the missing initial condition is given by $\frac{d^2f}{d\eta^2}(0) = -0.627555$.

\begin{figure}[!hbt]
	\centering
\psfrag{e}[l][]{$\eta$}  
\psfrag{f}[l][]{$f(\eta)$}  
\psfrag{df}[l][]{$\ds \frac{df}{d\eta}(\eta)$} 
\psfrag{ddf}[l][]{$\ds \frac{d^2f}{d\eta^2}(\eta)$} 
\includegraphics[width=13cm,height=13cm]{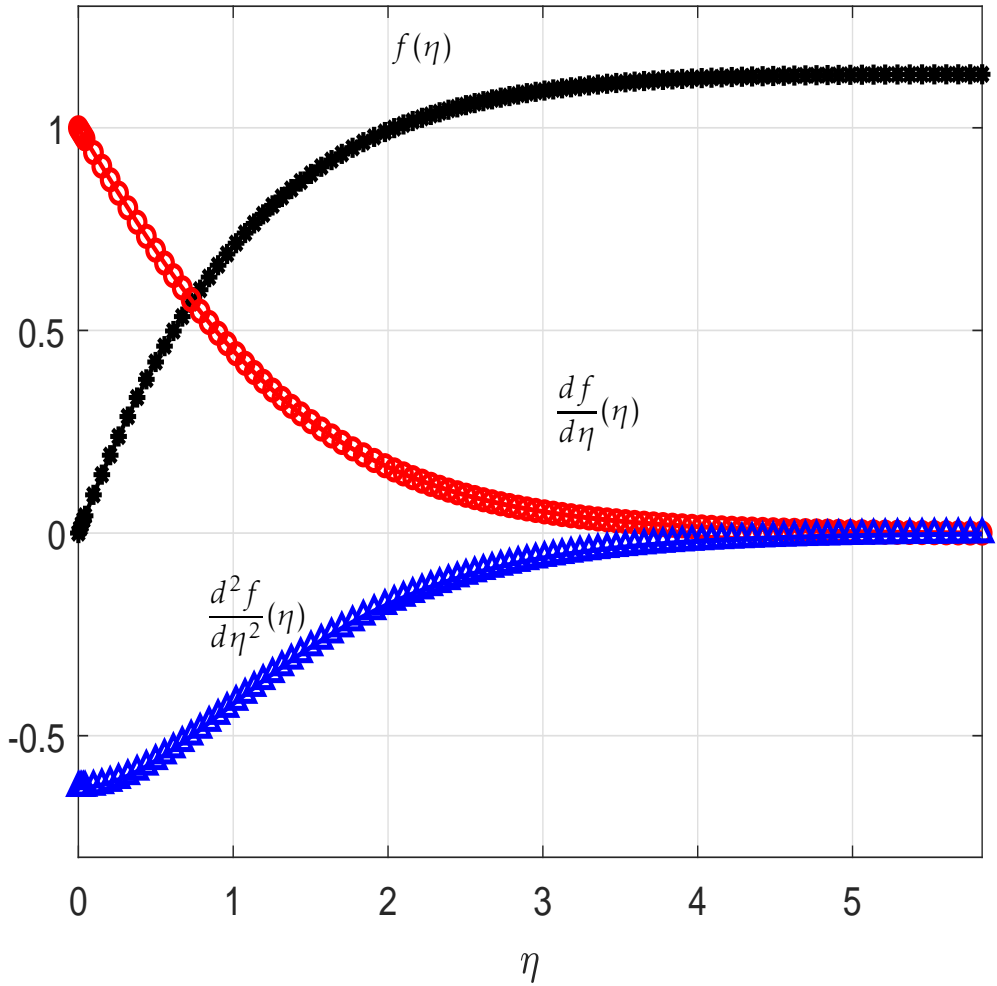}
\caption{Numerical results of the non-ITM for (\ref{eq:goveq}) with $\beta = 0$.}
	\label{fig:LiaoPopB0}
\end{figure}

Figure \ref{fig:LiaoPopB0} shows the solution of the Liao-Pop's problem in the particular case when we set $\beta = 0$.
As the considered mathematical model depends on a physical parameter we decided to list the obtained values of the missing initial condition versus this parameter.
In table \ref{tab:NITM:missingIC} we report the chosen parameter values and the related missing initial conditions $\frac{d^2f}{d\eta^2}(0)$.
\begin{table}[!hbt]
\caption{Numerical data and results.}
\vspace{.5cm}
\renewcommand\arraystretch{1.3}
	\centering
		\begin{tabular}{lr@{.}llr@{.}l}
\hline \\[-2.2ex]
{$\beta$} & \multicolumn{2}{c}%
{$ {\displaystyle \frac{d^2f}{d\eta^2}(0)}$} 
& {$\beta$} & \multicolumn{2}{c}%
{$ {\displaystyle \frac{d^2f}{d\eta^2}(0)}$}\\[1.5ex]
\hline
$-1$   & $-$0 & $258D-08$ & 0.1 & $-$0 & 672448 \\ 
$-0.9$ & $-$0 & 089014 & 0.2 & $-$0 & 712024 \\ 
$-0.8$ & $-$0 & 169015 & 0.3 & $-$0 & 754917 \\ 
$-0.7$ & $-$0 & 241805 & 0.4 & $-$0 & 793768 \\ 
$-0.6$ & $-$0 & 308699 & 0.5 & $-$0 & 831226 \\ 
$-0.5$ & $-$0 & 370678 & 0.6 & $-$0 & 867413 \\ 
$-0.4$ & $-$0 & 428499 & 0.7 & $-$0 & 902435 \\
$-0.3$ & $-$0 & 482755 & 0.8 & $-$0 & 936387 \\
$-0.2$ & $-$0 & 533922 & 0.9 & $-$0 & 969351 \\
$-0.1$ & $-$0 & 582389 & 1   & $-$1 & 001400 \\  
     0 & $-$0 & 628475 & & & \\ 
\hline			
		\end{tabular}
	\label{tab:NITM:missingIC}
\end{table}

\section{Concluding remarks} 
In this paper we have defined and applied an ITM for the numerical solution of the model studied by Liao and Pop \cite{Liao:2004:EAS}.
This method, that is based on the scaling invariance theory, solves a given BVP by solving a related sequence of IVPs and, therefore, it is an initial value method (like the shooting method).
We notice here that some exact solution were reported by Crane \cite{Crane:1970:FPS}.
For instance, when $\beta = 1$ he found the exact solution 
\[
f(\eta) = 1-exp(-\eta) \ , \qquad \frac{d^2f}{d\eta^2}(0) = 1 \ .
\]
Moreover, if $\beta = -1$ then he got the exact solution
\[
f(\eta) = \sqrt{2} \tanh \left(\eta/\sqrt{2}\right) \ , \qquad \frac{d^2f}{d\eta^2}(0) = 0 \ .
\]
Therefore, as it is easily seen the numerical results reported on table \ref{tab:NITM:missingIC} for these two special cases are just approximations.
We would like to remark here, how in the case $\beta = -1$ we are taking the ITM to its natural limit.
In fact, in this specific case we have to get a zero value rescaling a non-zero result. 
However, it is reassuring to see that the obtained value, namely $258D-08$, is very close to zero.
 
\vspace{1.5cm}

\noindent {\bf Acknowledgement.} {The research of this work was 
partially supported by the FABR grant of the University of Messina and by the GNCS of INDAM.}


\end{document}